\documentclass[a4paper,fleqn]{article}

\title{Some results in the exceptional set \\of \\Twin Prime Problem}
\author{Goldtwe  Anihc}
\date{August 1 ,2006}
\begin{document}
\maketitle
\newpage
 \begin{verbatim}
  Abstract

In the paper, there are new found methods to determine the range of
every exceptional element in  exceptional set, we can solve  Twin
primes problem and Goldbach Conjecture problem basically.


Key Word

exceptional set   Twin primes problem   Goldbach Conjecture problem


Address:   Gongnong Sicun  103--301
           Shanghai,P. R. of China

Email:  goldtwe@hotmail.com

 \end{verbatim}
\tableofcontents
\newpage
\section{{\LARGE\S1}}

  As the writer of "Sieve Methods" (i.e.[1]) indicates,
  for Twin primes problem and Goldbach Conjecture problem,
  they belong to the same problem.\\
  We can expand all conclusions (e.g. the result in [2],these results before the 1980s.) of Goldbach Conjecture problem
  to solve the problem of the twin prime, and vice-versa.\\

  The problem of  exceptional set exists in Goldbach Conjecture problem(cf.[2])
  after the three primes is proven. Therefor, the  problem of  exceptional set also exists in Twin primes problem,
  and their distribution of the exceptional element are same.\\
  We have the expression of  exceptional set in Goldbach Conjecture problem as
  follows.\\
     Let x be a larger positive integer.

\[\textrm{non-Goldbach number}:=\{n\not=p_1+p_2,2<p_1\le{n},2<p_2\le{n},2|n\}\]
\[\textrm{E}_g(x):=\{a:a = \textrm{ non-Goldbach number } n ,2<n\le{x}\}\]

and for the number of all element in  $\textrm{E}_g(x)$ set, we
express it as $\textrm{E}_g(x)$.\\

 In  $(x/2 ,  x]$, for all  even number n except $\textrm{E}_g(x)$ exceptional values,(cf.[2,3]) ,we have
\begin{equation}
 \textrm{D}_g(n)=2\textrm{C}(n)
 \frac{n}{\ln^{2}n}+O(\frac{x(\ln\ln x)^{3}}{\ln ^{3}x})
\end{equation}
where
\[\textrm{C}(n)=\prod_{p>2}(1-\frac{1}{{(p-1)}^{2}})\prod_{p>2}^{p|n}\frac{p-1}{p-2}\]
and $\textrm{D}_g(n)$,cf.(10)
\begin{equation}
\textrm{E}_g(x)\ll\frac{x}{\ln^{A}x}
\end{equation}
where, $A$ be an any give positive, and the contained constant
depending on $A$. \\Here, In Twin primes problem,

\[\textrm{C}(n)=\prod_{p>2}(1-\frac{1}{{(p-1)}^{2}})\]

For the  exceptional set of Twin primes problem, to the  sieve
function eye,we have
\[\textrm{A}_g:=\{a:a = n-p ,2<p\le{n},2|n\}\]

where $\textrm{A}_g$  which is in the  sieve function
$\textrm{S}_g(\textrm{A}_g, P , z)$ of Goldbach Conjecture problem
be a set.

\[\textrm{A}_t:=\{a:a = p+2 ,2<p\le{n},2|n\}\]

where $\textrm{A}_t$ which is in the  sieve function
$\textrm{S}_t(\textrm{A}_t, P , z)$ of Twin primes problem be a set.

    As for $\textrm{E}_g(x)$,an  even number $n$ is an exceptional
element when all $n - p$ are not the prime in $\textrm{A}_g$.

   To $\textrm{E}_t(x)$ of Twin primes problem, it is the same with
$\textrm{E}_g(x)$ of Goldbach Conjecture problem. i.e. An  even
number $n$ is an exceptional element  when all $p + 2$ are not the
prime in $\textrm{A}_t$. This conclusion developed from the
$\textrm{E}_g(x)$ of Goldbach Conjecture problem. Their distribution
of the exceptional element are same.
   On this, we have
\begin{equation}
 \textrm{D}_t(n)=2\prod_{p>2}(1-\frac{1}{(p-1)^{2}})
 \frac{n}{\ln^{2}n}+O(\frac{x(\ln\ln x)^{3}}{\ln^{3}x})
\end{equation}
where $\textrm{D}_t(n)$, cf.(11).
\begin{equation}
 \textrm{E}_t(x)\ll\frac{x}{\ln^{A}x}
\end{equation}
where, $A$ be an any give positive, and the contained constant
depending on $A$.

\section{{\LARGE\S2}}

  We quote the increment to the  sieve function of Twin primes problem. And the fundamental
properties of  sieve function as follows.
\begin{equation}
 \textrm{(i)}S(A,P,z_1)\leq S(A,P,z_2)\qquad,\qquad(z_1\geq z_2)
\end{equation}
\begin{equation}
 \textrm{(ii)} 0\leq S(A,P,z)\leq |A|
\end{equation}
\begin{equation}
 \textrm{(iii)} S(A+\Delta A,P,z)= S(A,P,z)+S(\Delta A,P,z)
\end{equation}
Here, $|A|$ be a number of all element in set $A$.
      $\Delta A$ be the non-empty subset of  set \mbox{$(A+\Delta A)$},and we also call it
after an increment of set  $A$. Its number is $|\Delta A|$.
  To $A$ and $\Delta A$, their cap is an empty set $\phi$.i.e.

\[A\cap \Delta A=\phi\]
  Then
\begin{eqnarray}
 S(A+\Delta A,P,z_1) & = & S(A,P,z_1)+S(\Delta A,P,z_1)\nonumber \\
                    & < & S(A,P,z_2)+|\Delta A|
\end{eqnarray}
Where  $z_1> z_2 $.\\
  So,for three expressions in (8), if we know two expressions of they, then another of they
has the upper or lower bound.\\
  When we apply (8), (3), (4) and Drawer principle to the  exceptional set
of Twin primes problem, for \mbox{$2(\textrm{E}_t(x) +1)$} natural
numbers in a closed interval of the natural number,
\mbox{$\{2|M,[M,M+2\textrm{E}_t(x)],(x/2<M, M+2\textrm{E}_t(x)\leq
x]\}$},if an even number which is one of two endpoints is an
exceptional element,then the non-exceptional element which is an
even number in this closed interval must exists. And the difference
between this two  even numbers less than or equal to
$2\textrm{E}_t(x)$. Besides, for $\S3 $,the exceptional element is
the right endpoint.
  On this
\[|\Delta A|\leq2\textrm{E}_t(x)\]
and $S( A, P, z_1)$ or $S( A, P, z_2)$ which one of they in (8) be
the $\textrm{D}_t(n)$ of (3).
  At this very moment, $A>4$ in (4),we have

\[\frac{|\Delta A|}{\textrm{D}_t(n)}\longrightarrow 0\qquad,\qquad x\rightarrow \infty\]

  So, we must know that an exceptional element of Twin primes problem
which is the another of they has the upper or lower bound.
  As for an exceptional element of Goldbach Conjecture problem, by the $\S1 $,
it also has the upper or lower bound.

\section{{\LARGE\S3}}

  When  even number $n_1>n_2$, if the expressions of $\textrm{D}_t(n_1)-\textrm{D}_t(n_2)$
as (15), (16) and (17), then we can determine the range of every
exceptional element in exceptional set of Twin primes problem.
  So, for $\alpha$ and $\beta$ sequence of number as follows.\\
\begin{verbatim}
          i        1   2   3  .......  n-1   n
        \alpha     1   2   3  .......  n-1   n ,  2|n
        \beta_g    n-1 n-2 n-3 .......   1
        \beta_t     3   4   5  .......  n+1  n+2
\end{verbatim}

  First, by Eratosthenes' sieve, we have all primes in $\alpha$ and $\beta$ sequence of number.
  Second, we map the process that sieve $\beta$ sequence of number into $\alpha$ sequence
of number.
  The manner of mapping be \mbox{$n_{\alpha i} + n_{\beta_g  i}=n$}
or \mbox{$n_{\beta_t i}-n_{\alpha i}= 2$}. And we also sieve out
these number which was mapped in a sequence of number.
  Then, we  have (10) and (11).
       Let  $P$ is  the prime set in $[2 , n]$,  its number $|P|$ is $\pi(n)$. i.e.
 \[P : = \{ p : 1< p \leq n \}\qquad, \qquad|P| = \pi(n)\]
  For $\alpha$ sequence of number, we have the prime theorem as follows.
\begin{equation}
 \pi(n)=n+\sum_{d|\prod_{p\in P} p}\mu(d)[\frac{n}{d}]+O(\pi(\sqrt{n}))
\end{equation}
After  the mapping, for Goldbach Conjecture problem

\[\textrm{D}_g:=\{ a\in P :  a = n - p,  2<p\leq n, 2|n \}\]

 \[\textrm{D}_g(n)=|\textrm{D}_g| \]
 (cf.$P$ \quad and \quad $|P|$)
 i.e.
\begin{equation}
 \textrm{D}_g(n)=\pi(n)+\sum_{d|\prod_{p\in P}p}^{2|n}\mu(d)\pi(n,d,{n}\!\!\!\!\pmod{d})+O(\pi(\sqrt{n}))
\end{equation}

  By the same manner,  for Twin primes problem
\[\textrm{D}_t:=\{ a\in P :  a = p+2,  2<p\leq n, 2|n \}\]
\[\textrm{D}_t(n)=|\textrm{D}_t|\]
\textrm{(cf.$P$\quad and\quad $|P|$)} i.e.
\begin{equation}
 \textrm{D}_t(n)=\pi(n)+\sum_{d|\prod_{p\in P}p}^{2|n}\mu(d)\pi(n,d,d-2) +O(\pi(\sqrt{n}))
\end{equation}
  It may be seen ,when \(n\rightarrow  \infty \) ,we have the following

 \begin{equation}
  \textrm{D}_g(n)\longrightarrow\pi(n)+\sum_{d|\prod_{p\in P}p}^{2|n}\mu(d)\pi(n,d,{n}\!\!\!\!\pmod{d})
 \end{equation}
\begin{equation}
 \textrm{D}_t(n)\longrightarrow\pi(n)+\sum_{d|\prod_{p\in P}p}^{2|n}\mu(d)\pi(n,d,d-2)
\end{equation}

where $d-2$ of the second term does not change with  $n$ ,and this
second term be the element number of  a union of set. i.e.

\begin{equation}
 -\sum_{d|\prod_{p\in P}p}^{2|n}\mu(d)\pi(n,d,d-2)=|\bigcup_{2<p\leq\sqrt{n}}\{ a\in P: a \equiv {p-2}\!\!\!\!\pmod{p} \}|
\end{equation}
 (cf. $P$\quad  and \quad  $|P|$).

  So, when   even number $n_1 > n_2$ , Let

\[P_1:= \{ p: 2<p\leq n_1\}\]
\[P_2:= \{ p: 2<p\leq n_2\}\]
\[\Delta P:= \{ p: n_2<p\leq n_1\}\]

Then

\begin{eqnarray}
 \textrm{D}_t(n_1)-\textrm{D}_t(n_2) &  \longrightarrow &  \pi(n_1)-\pi(n_2)-{}
                                                                       \nonumber\\*
                                        &  &  {}-|\bigcup_{2<p\leq\sqrt{n_1}}\{a\in P_1:a \equiv {p-2}\!\!\!\!\pmod{p}\}|+{}
                                                                                                             \nonumber\\*
                                       &  &  {}+|\bigcup_{2<p\leq\sqrt{n_2}}\{ a\in P_2: a\equiv {p-2}\!\!\!\!\pmod{p}\}|
\end{eqnarray}

(cf.(13))

  Suppose,

\[\textrm{set1}:=\{a\in P_2 :a\equiv {p-2}\!\!\!\!\pmod{p},\sqrt{n_2}<p\leq \sqrt{n_1}\}\]

\[\textrm{set2}:=\{a\in \Delta P :a\notin \Delta P\cap
          \{b\in \Delta P:b\equiv {p-2}\!\!\!\!\pmod{p},2<p\leq \sqrt{n_1}\}\}\]
  Then
\begin{eqnarray}
 \textrm{D}_t(n_1)-\textrm{D}_t(n_2) &  \rightarrow &  -|\textrm{set1}|+|\textrm{set2}|
                                                                                 \nonumber\\*
                                       &   < & |\textrm{set1}|+|\textrm{set2}|
\end{eqnarray}

(cf. $P$ and $|P|$),where
\[\textrm{set1}\in\bigcup_{2<p\leq\sqrt{n_1}}\{a\in P_1:a\equiv {p-2}\!\!\!\!\pmod{p}\}\]

  Now we evaluate the value of $|\textrm{set1}|$ and $|\textrm{set2}|$.
\[0\leq \textrm{set2}<|\Delta P|<n_1-n_2 \qquad ,
                                                   \qquad (cf.(6),P \quad and \quad |P|)\]
\[0\leq \textrm{set1} \  < \  \frac{n_2}{p_{min}}(\sqrt{n_1}-\sqrt{n_2}){}
                                                                  \nonumber\\*
                           \  \  {} <[(\sqrt{n_1n_2})-n_2]{}
                                                                  \nonumber\\*
                           \  \  {} < (n_1-n_2)\]
where $p_{min}\rightarrow\sqrt{n_2}$,$\quad$
cf.$\sqrt{n_2}<p\leq\sqrt{n_1}$

  So,
\begin{equation}
 0\leq\textrm{D}_t(n_1)-\textrm{D}_t(n_2)<O(n_1-n_2)
\end{equation}

  When  $A>4$ in (4), every exceptional element in Twin primes problem
are evaluated as the ending part of $\S2 $. Thus we solve Twin
primes
problem basically.\\
  By the same manner, every exceptional element in Goldbach Conjecture problem are also
evaluated. Thus we also solve Goldbach Conjecture problem basically.
(end)

\section{References}
 \begin{verbatim}
[1]H.Halberstam, H.E.Richert, Sieve Methods, Academic Press,London,
                              1974.

[2] T.Estermann, Proof that almost all even positive integers are
                 sun of two primes,
                 Proc. London Math. Soc.(1).44(1938), 307--314.

[3] C.D.Pan,C.B.Pan ,Goldbach Conjecture, Science Press, Beijing,
                     1992;(in Chinese) 1984.

(end)
 \end{verbatim}
\end{document}